\documentclass[11pt]{article}
\usepackage{geometry}
\newcommand{\cupall}{\pmb{\pmb{\cup}}}

\usepackage[T2A]{fontenc}
\usepackage{amsmath,amsthm,amsfonts,amssymb,epsfig,color}
\usepackage[greek,english]{babel}
\usepackage[utf8]{inputenc} 
\usepackage{graphicx}  
\usepackage{bm}  

\usepackage{enumerate}

\usepackage[bookmarks,breaklinks]{hyperref}  

\newcommand{\tw}{{\mathbf{tw}}}

\newcommand{\eg}{{\mathbf{eg}}}
\newcommand{\an}{{\mathbf{an}}}
\newcommand{\n}{{\mathbf{n}}}

\theoremstyle{remark}

\theoremstyle{plain}
\newtheorem{corollary}{Corollary}
\newtheorem{theorem}{Theorem}
\newtheorem{lemma}{Lemma}
\newtheorem{observation}{Observation}
\newtheorem{definition}{Definition}
\newtheorem{proposition}{Proposition}

\begin{document}
\title{Optimizing  the Graph Minors Weak Structure Theorem}
\author{Archontia C. Giannopoulou\thanks{Department of Mathematics, National and Kapodistrian University of Athens, Greece. }~\thanks{Supported by a grant of the  Special Account for Research Grants of the National and Kapodistrian University of Athens (project code: 70/4/10311).} \and Dimitrios M. Thilikos$^*$}

\date{}
\maketitle
\begin{abstract}
\noindent 
One of the major results of  [{\sl N. Robertson and P.~D. Seymour. Graph minors. {XIII}. {T}he disjoint paths problem.} {\em J. Combin. Theory Ser. B}, {\sl 63(1):65--110, 1995}], also known as {\sl the weak structure theorem}, revealed  the 
local  structure of graphs excluding some graph as a minor: each such graph $G$ either has small treewidth or contains the subdivision of a wall 
that can be arranged ``bidimensionally'' inside $G$, given that 
some small set of vertices are removed. We prove  an optimized version of 
that theorem where (i)
the relation between the treewidth of the graph and the height of the wall is linear (thus best possible) and (ii)
the number of vertices to be removed is minimized.\\

\end{abstract}

\section{Introduction}

The Graph Minors series of Robertson and Seymour appeared to be a rich source 
of structural results in graph theory with multiple applications in algorithms. One of the most 
celebrated outcomes of this project was the existence of an $O(n^{3})$ step algorithm 
for solving problems such as the {\sc Disjoint Path} and the {\sc Minor Containment}.
A basic ingredient of this algorithm is a theorem, proved in paper XIII of the series~\cite{RobertsonS-XIII},
revealing the local structure of graphs excluding some graph as a minor. This result, now called 
the {\em weak structure theorem}, asserts that there is a function $f: \Bbb{N}\times \Bbb{N}\rightarrow \Bbb{N}$ such that for every integer $k$, every $h$-vertex graph $H$, and every graph $G$, one of the following holds:
\begin{itemize}
\item[1.] $G$ contains $H$ as a minor, 
\item[2.] $G$ has treewidth at most $f(k,h)$, or  
\item[3.] $G$ contains a set $X$ of at most ${h \choose 2}$ vertices (called {\em apices}) such that $G \setminus X$
contains as a subgraph the subdivision $W$ of a  wall of height $k$ that is arranged inside $G$ in a ``flat'' manner ({\sl flatness condition}).
\end{itemize}
To make the above statement precise we need to clarify the flatness condition in the third statement above. We postpone this complicated task until Section~\ref{sec:fsw} and instead, we roughly visualize $W$ in a way that the part of $G\setminus X$ that is 
located inside the perimeter $P$ of $W$ can be seen as a set of graphs attached on a plane region where each of these graphs has bounded treewidth and its boundary with the other graphs is bounded by 3.

The algorithmic applications of the weak structure  theorem reside in  the fact that the graph 
inside $P$ can be seen as a {\sl bidimensional structure} where, for several combinatorial problems, a solution certificate can be revised so that it avoids the middle vertex 
of the subdivided wall $W$. 
This is known as the {\em irrelevant vertex technique} and can be seen as a reduction of  an instance of a problem to an equivalent one where this ``irrelevant vertex'' has been deleted.
The application of this technique has now gone much further than its original use in 
the Graph Minors series and has evolved to a standard tool in algorithmic graph minors theory (see~\cite{DawarGK07,DawarK09,GolovachKPT09,KawarabayashiK08,Kawarabayashi:2010cs,Kobayashi:2009jt} for applications of this technique).

In this paper we prove an optimized version of the weak structure theorem.
Our improvement is twofold: first, the function $f$ is now linear on $k$ and second,
the number of the apices is bounded by $h-5$.
Both our optimizations are 
optimal as indicated by the graph $J$ obtained by taking 
a $(k\times k)$-grid (for $k\geq 3$) and making all its vertices adjacent with a copy of $K_{h-5}$.
Indeed, it is easy to verify that $J$ excludes $H=K_{h}$ as a minor, its treewidth is $k+h-5$ and becomes planar 
(here, this is equivalent to the flatness condition)
after the removal of exactly $h-5$ vertices.

Our proof deviates significantly from the one in~\cite{RobertsonS-XIII}. It
builds on the (strong) {\sl structure theorem} of the Graph Minors that was proven in paper XVI of the series~\cite{RobertsonS-XVI}. This theorem reveals the global structure of a graph 
without a $K_{h}$ as a minor and asserts  that each such graph can be obtained 
by gluing together graphs that can ``almost be embedded'' in a surface where 
$K_{h}$ cannot be embedded (see~\ref{sed:prel} for the exact statement).
The proof exploits this structural result  and combines it  with the fact, proved in~\cite{FominGT09con}, that apex-free
``almost embedded graphs'' without a $(k\times k)$-grid have treewidth  $O(k)$.

The organization of the paper is the following. In Section~\ref{sed:prel} we give the 
definitions of all the tools that we are going to use in our proof, including the 
 Graph Minors structure theorem. The definition of the flatness condition is given
 in Section~\ref{sec:fsw}, together with the statement of our main result. Some lemmata concerning the invariance of the flatness property under certain local transformations are given in Section~\ref{fil} and further definitions and results 
concerning apex vertices are given in Section~\ref{sec:pytr}. Finally, the proof 
of our main result is presented in Section~\ref{sec:proof}.

%
%
%
%
%
%
%

\section{Preliminaries}
\label{sed:prel}

Let $n\in \mathbb{N}$. We denote by $[n]$ the set $\{1,2,\ldots, n\}$. Moreover, for every $k\leq n$, if $S$ is a set such that $|S|=n$ we say that a set $S'\subseteq S$ is a $k$-subset of $S$ if $|S'|=k$.

A {\em graph} $G$ is a pair $(V,E)$ where $V$ is a finite set, called the {\em vertex set} and denoted by $V(G)$, and $E$ is a set of 2-subsets of $V$, called the {\em edge set} and denoted by $E(G)$. We denote by $\n(G)$ the number of vertices of $G$, i.e. $\n(G)=|V(G)|$.
If we allow $E$ to be a subset of $\mathcal{P}(V)$ then we call the pair $H=(V(H),E(H))$ a {\em hypergraph}. The {\em incidence graph} of a hypergraph $H$ is the bipartite graph $I(H)$ on the vertex set $V(H)\cup E(H)$ where $v\in V(H)$ is adjacent to $e\in E(H)$ if and only if $v\in e$, i.e. $v$ is incident to $e$ in $H$. We say that a hypergraph $H$ is {\em planar} if its incidence graph is planar. Unless otherwise stated, we consider finite undirected graphs without loops or multiple
edges. 

Let $G$ be a graph. For a vertex $v$, we denote by $N_G(v)$ its
\emph{(open) neighborhood}, i.e. the set of vertices which are
adjacent to $v$. The \emph{closed neighborhood} $N_G[v]$ of $v$ is the set
$N_G(v)\cup\{v\}$. For $U\subseteq V(G)$, we
define  $N_G[U]=\bigcup_{v\in U}N_G[v]$. We may omit the index if the
graph under consideration is clear from the context. If $U\subseteq
V(G)$ (resp. $u\in V(G)$ or $E\subseteq E(G)$ or $e\in E(G)$) then
$G-U$ (resp. $G-u$ or $G-E$ or $G-e$) is the graph obtained from $G$
by the removal of vertices of $U$ (resp. of vertex $u$ or edges of
$E$ or of the edge $e$). We say that a graph $H$ is a subgraph of a graph $G$, denoted by $H\subseteq G$, if $H$ can be obtained from $G$ after deleting edges and vertices. Let ${\cal C}$ be a class of graphs and $S$ be a set of vertices. We denote by $\cupall {\cal C}$ the graph $\cup_{G\in \mathcal{C}}G$ and by ${\cal C}\setminus S=\{G\setminus S\mid G\in {\cal C}\}$.

The complete graph on $n$ vertices is denoted by $K_{n}$. Moreover, if $S$ is a finite set, we denote by $K[S]$ the complete graph with vertex set $S$.
Let $G$ be a graph such that $K_{3}\subseteq G$ and $x,y,z$ be the vertices of $K_{3}$. The {\em $\Delta Y$-transformation of $K_{3}$ in $G$} is the following: We remove the edges $\{x,y\},\{y,z\},\{x,z\}$, add a new vertex $w$ and then add the edges $\{x,w\},\{y,w\},\{z,w\}$. 

Let $G$ be a graph. We say that $G$ is an apex graph if there exists a vertex $v\in V(G)$ such that $G\setminus v$ is planar. Moreover, we say that $G$ is an $\alpha_{G}$-apex graph if there exists an $S\subseteq V(G)$ such that $|S|\leq \alpha_{G}$ and $G\setminus S$ is planar. We denote by $\an(G)$, the minimum $k\in \mathbb{N}$ such that $G$ is a $k$-apex graph, i.e. $\an(G)=\min\{k\in \mathbb{N}\mid \exists S\subseteq V(G):(|S|\leq k\land G\setminus S \text{ is planar})\}$. Clearly, $\mathcal{G}=\{G\mid\an(G)=1\}$ is the class of the apex graphs.

\begin{observation}\label{verwithungrtw}
Let $T$ be a tree, $k\in \mathbb{N}$ and $w:V(T)\rightarrow \mathbb{N}$ such that there exists at least one vertex $v\in V(T)$ with $w(v)\geq k$. There exists a vertex $u\in V(T)$ such that at most one of the connected components of $(G\setminus u)$ contains a vertex $u'$ with $w(u')> k$.
\end{observation}

\begin{proof}
Let $Y=\{v\in V(T)\mid w(v)\geq k\}$. Pick a vertex $r$ of $T$ and let 
$v$ be a vertex of $Y$ with maximum distance away from $r$. It is easy to verify that the lemma holds for $v$.
\end{proof}

\paragraph{\bf Surfaces.}

A \emph{surface} $\Sigma$ is a compact 2-manifold without boundary
(we always consider connected surfaces).
Whenever we refer to a {\em
$\Sigma$-embedded graph} $G$ we consider a  2-cell embedding of
$G$ in $\Sigma$. To simplify notations, we do not distinguish
between a vertex of $G$ and the point of $\Sigma$ used in the
drawing to represent the vertex or between an edge and the line
representing it.  We also consider a graph $G$ embedded in
$\Sigma$ as the union of the points corresponding to its vertices
and edges. That way, a subgraph $H$ of $G$ can be seen as a graph
$H$, where $H\subseteq G$.
Recall that $\Delta \subseteq \Sigma$ is
an open  (resp. closed)  disc if it is homeomorphic to
$\{(x,y):x^2 +y^2< 1\}$ (resp. $\{(x,y):x^2 +y^2\leq 1\}$).
The {\em Euler genus} of a non-orientable surface $\Sigma$
is equal to the non-orientable genus
$\tilde{g}(\Sigma)$ (or the crosscap number).
The {\em Euler genus}  of an orientable   surface
$\Sigma$ is $2{g}(\Sigma)$, where ${g}(\Sigma)$ is  the orientable genus
of $\Sigma$. We refer to the book of Mohar and Thomassen \cite{MoharT01} for
more details  on graphs embeddings.
The {\em Euler genus} of a graph $G$ (denoted by $\eg(G)$) is the minimum integer $\gamma$ such that $G$ can be embedded on a surface of the Euler genus $\gamma$.

\paragraph{\bf Contractions and minors.}

Given an edge  $e=\{x,y\}$ of a graph $G$, the graph  $G/e$ is
obtained from  $G$ by contracting the edge $e$, i.e.
the endpoints $x$  and $y$ are replaced by a new vertex $v_{xy}$
which  is  adjacent to the old neighbors of $x$ and $y$ (except
$x$ and $y$). 
A graph $H$ obtained by a sequence of
edge-contractions is said to be a \emph{contraction} of $G$. An alternative, and more useful for our purposes, definition of a contraction is the following.

Let $G$ and $H$ be graphs and let $\phi: V(G)\rightarrow V(H)$ be
a surjective mapping  such that
\begin{itemize}
\item[${\bf 1.}$]  for every vertex $ v\in V(H)$, its codomain $\phi^{-1}(v)$ induces connected graph
$G[\phi^{-1}(v)]$; 
\item[${\bf 2.}$]\label{phicon} for every edge $ \{v,u\}\in E(H)$, the graph
$ G[\phi^{-1}(v)\cup \phi^{-1}(u)]$ is
connected;
\item[${\bf 3.}$] for every $ \{v,u\}\in E(G)$,  either $\phi(v)=\phi(u)$,
or $\{\phi(v),\phi(u)\}\in E(H)$.
\end{itemize}
We, then, say that {\em $H$ is a contraction of $G$ via $\phi$}
and  denote it as $H\leq_{c}^{\phi}G$. Let us observe  that
{$H$ is a contraction of $G$}  if
$H\leq_{c}^{\phi} G$ for some  $\phi: V(G)\rightarrow V(H)$.
In this case we  simply write $H\leq_{c} G$.
If $H\leq_{c}^{\phi} G$ and $v\in V(H)$, then
we call the codomain $\phi^{-1}(v)$ the {\em model of $v$} in $G$.

Let $G$ be a graph embedded in some surface $\Sigma$ 
and let  $H$ be a contraction of $G$ via
function $\phi$. We say  that $H$ is a \emph{surface
contraction} of $G$ if for each vertex $v\in V(H)$,
$G[\phi^{-1}(v)]$ is embedded in some open disk in $\Sigma$.

Let $G_{0}$ be a graph embedded in some surface $\Sigma$ of Euler genus $\gamma$ and
let $G^{+}$ be another graph that might share common vertices with $G_{0}$. We set
$G=G_{0}\cup G^{+}$. Let also $H$ be some graph and let $v\in V(H)$.
We say that $G$ {\em contains a graph $H$ as a $v$-smooth contraction}
if $H\leq_{c}^{\phi}G$ for some $\phi:V(G)\rightarrow V(H)$ and
there exists a closed  disk $D$ in
$\Sigma$ such that all the vertices of $G$ that are outside $D$ are
exactly the model of $v$, i.e. $\phi^{-1}(v)= V(G)\setminus (V(G)\cap D)$.

A graph $H$ is a {\em minor} of a graph $G$, denoted by $H\leq_{m} G$, if $H$ is the contraction of some subgraph
of $G$. If we restrict the contraction to edges whose one of the incident vertices has degree exactly two, also called {\em disolving} that vertex, then we say that $H$ is a {\em topological minor} of $G$ and we denote it by $H\leq_{tm} G$.  Moreover, we say that a graph $G$ is a {\em subdivision} of a graph $H$, if $H$ can be obtained from $G$ by disolving vertices.

We say that a graph $G$ is {\em $H$-minor-free} when it does not
contain $H$ as a minor. We also say that a graph class ${\cal G}$
is {\em $H$-minor-free} (or, excludes $H$ as a minor)  when
all its members are $H$-minor-free.

\paragraph{\bf Graph Minors structure theorem.}

The proof of our results is using  the Excluded Minor Theorem from
the Graph Minor theory. Before we state it, we need some
definitions.

\begin{definition}[{\rm{\em $h$-nearly embeddable graphs}}]
Let $\Sigma$ be a surface and $h>0$ be an integer.  A graph
$G$ is $h$-nearly embeddable in $\Sigma$ if there is a set of  vertices 
$X\subseteq V(G)$   (called apices)  of size at most $h$ such that graph $G- X$ is the union of 
 subgraphs $G_0,\dots,G_h$ 
with the following properties
\begin{enumerate}[i)]
\item There is a set of 
 cycles $C_1, \dots,C_h$ in $\Sigma$  such that the cycles $C_i$ are the borders of open pairwise disjoint discs $\Delta_{i}$ in $\Sigma$;
\item $G_0$ has an  embedding in $\Sigma$ in such a way that $G_{0}\cap \bigcup_{i=1,\ldots,h}\Delta_{i}=\emptyset$;
\item graphs $G_1,\dots,G_h$ (called vortices) are pairwise disjoint and  for $1\leq i \leq h$, $V(G_0)\cap V(G_i)\subset C_i$;
\item for $1\leq i \leq h$, let $U_i:= \{u_{1}^{i},\dots,u_{m_{i}}^{i}\} $ be the vertices of $V(G_0) \cap V(G_i)\subset C_i$ appearing in an order obtained by  clockwise traversing of $C_i$, we call vertices of  $U_i $ {\em bases} of $G_{i}$. Then $G_i$ has a
path decomposition $\mathcal{B}_{i}=(B_{j}^{i})_{1\leq j \leq m_i}$, of width at most $h$ such that
 for $1 \leq j \leq m_i$, we have $u_j^{i} \in B_{j}^{i}$.
\end{enumerate}
\end{definition}

%

\noindent A \emph{tree decomposition} of a graph $G$ is a pair $(\mathcal{X},T)$ where $T$
is a tree and ${\cal X}=\{X_{i} \mid i\in V(T)\}$ is a collection of subsets
of $V(G)$ such that:
\begin{itemize}
\item[{\bf 1}.]  $\bigcup_{i \in V(T)} X_{i} = V(G)$;
\item[{\bf 2}.] for each edge $\{x,y\} \in E(G)$, $\{x,y\}\subseteq X_i$ for some
$i\in V(T)$, and
\item[{\bf 3}.]\label{indconsub} for each $x\in V(G)$ the set $\{ i \mid x \in X_{i} \}$
induces a connected subtree of $T$.
\end{itemize}
The \emph{width} of a tree decomposition $(\{ X_{i} \mid i \in V(T) \},T)$
 is $\max_{i \in V(T)}\,\{|X_{i}| - 1\}$. The \emph{treewidth} of a
graph $G$ is the minimum width over all tree decompositions of $G$.
Furthermore, we call the subsets $X_{i}$ {\em bags} of the decomposition and for every $X_{i},i \in V(T)$ we denote by $\overline{X}_{i}$ its closure, i.e. $\displaystyle \overline{X}_{i}$ is the graph $\displaystyle G\left[X_{i}\right]\cup\left(\cup_{j\in N_{T}(i)}K\left[X_{i}\cap X_{j}\right]\right)$.

\begin{observation}
\label{treew:bags}
Let $G$ be a graph and $(\mathcal{X},T)$ be a tree decomposition of $G$. Then there exists an $X\in \mathcal{X}$ such that $\tw(\overline{X})\geq \tw(G)$.
\end{observation}

We also need the simple following result.

\begin{lemma}
\label{apex} If $G$ is a graph and $X\subseteq V(G)$, then $\tw(G-
X)\geq \tw(G)-|X|$.
\end{lemma}

We say that a tree decomposition $(\mathcal{X},T)$ of a graph $G$ is {\em small} if for every $i,j\in V(T)$, $i\neq j$, $X_{i} \nsubseteq X_{j}$.

A simple proof of the following lemma can be found in~\cite{FlumG06par}.

\begin{lemma}\label{fincomptd}
\begin{enumerate}
\item Let $G$ be a graph and $(\mathcal{X},T)$ be a small tree decomposition of $G$. Then $|V(T)|\leq |V(G)|$.
\item Every graph $G$ has a small tree decomposition of width $\tw(G)$. 
\end{enumerate}
\end{lemma}
The following proposition is known as the Graph Minors structure  theorem  \cite{RobertsonS-XVI}.
\begin{proposition}
\label{structurethm} There exists a computable function $f:\mathbb{N}\rightarrow\mathbb{N}$ such that, for every non-planar graph $H$ with $h$ vertices and every graph $G$
excluding $H$ as a minor there exists a tree decomposition $(\mathcal{G}=\{G_{i}\mid i\in V(T)\},T)$ where for every $i\in V(T)$, $\overline{G}_{i}$ is an $f(h)$-nearly embeddable graph in a surface $\Sigma$ on which $H$ cannot be embedded.
\end{proposition}




\section{Statement of the main result}
\label{sec:fsw}

Let $k$ and $r$ be positive integers where $k,r\geq 2$. The
\emph{$(k\times r)$-grid} is the Cartesian product of two paths of
lengths $k-1$ and $r-1$ respectively. A vertex of a  $(k\times
r)$-grid is a {\em corner} if it has degree $2$. Thus each
$(k\times r)$-grid has 4 corners. A vertex of a  $(k\times
r)$-grid is called {\em internal} if it has degree 4, otherwise it
is called {\em external}.

We define $\Gamma_{k}$ as the following
(unique, up to isomorphism)
triangulation  of the $(k\times k)$-grid.
Let $\Gamma$ be a plane embedding of the $(k\times k)$-grid such that
all external vertices are on the boundary of the external face.
We triangulate internal faces of the $(k\times k)$-grid
such
that, in the obtained
graph, all the internal vertices have degree $6$ 
and all non-corner external vertices have degree 4. The construction of 
$\Gamma_{k}$ is completed if we connect one corner of degree two 
with all vertices of the external face (we call this corner {\em loaded}).
We also use notation $\Gamma_{k}^{*}$ for the graph obtained from $\Gamma_{k}$ if we remove all edges incident to its loaded vertex that do not exist in its underlying grid. 

We define the $(k,l)$-pyramid to be the graph obtained if we take the disjoint union of a $(k\times k)$-grid and a clique $K_{l}$ and then add all  edges between the vertices of the clique and the vertices of the grid. We denote the $(k,
l)$-pyramid by $\Pi_{k,l}$.

\paragraph{Walls.}
A \emph{wall of height $k$}, $k\geq 1$, is the graph obtained from a
$((k+1)\times (2\cdot k+2))$-grid with vertices $(x,y)$,
$x\in\{1,\dots,2\cdot k+4\}$, $y\in\{1,\dots,k+1\}$, after the removal of the
``vertical'' edges $\{(x,y),(x,y+1)\}$ for odd $x+y$, and then the removal of
all vertices of degree 1. 
We denote such a wall by $W_{k}$.

The  {\em corners} of the wall $W_{k}$ are the vertices $c_{1}=(1,1)$, $c_{2}=(2\cdot k+1,0)$, $c_{3}=(2\cdot k + 1 + (k+1\mod 2),k+1)$ and $c_{4}=(1+(k+1\mod 2),k+1)$. We let $C=\{c_{1},c_{2},c_{3},c_{4}\}$ and we call the pairs $\{c_{1},c_{3}\}$ and $\{c_{2},c_{4}\}$ {\em anti-diametrical}.

 A {\em subdivided wall $W$ of height $k$} is a wall obtained from $W_{k}$ after replacing some of its edges by paths without
common internal vertices. We call the resulting graph $W$ a {\em subdivision} of $W_{k}$ and the vertices that appear in the wall after the replacement {\em subdivision} vertices.

The non-subdivision vertices of $W$ are called  {\em original} vertices. 
The {\em perimeter} $P$ of a subdivided wall is the cycle defined by its boundary.

\begin{figure}[h]
  \begin{center}
\includegraphics[width=10.18cm]{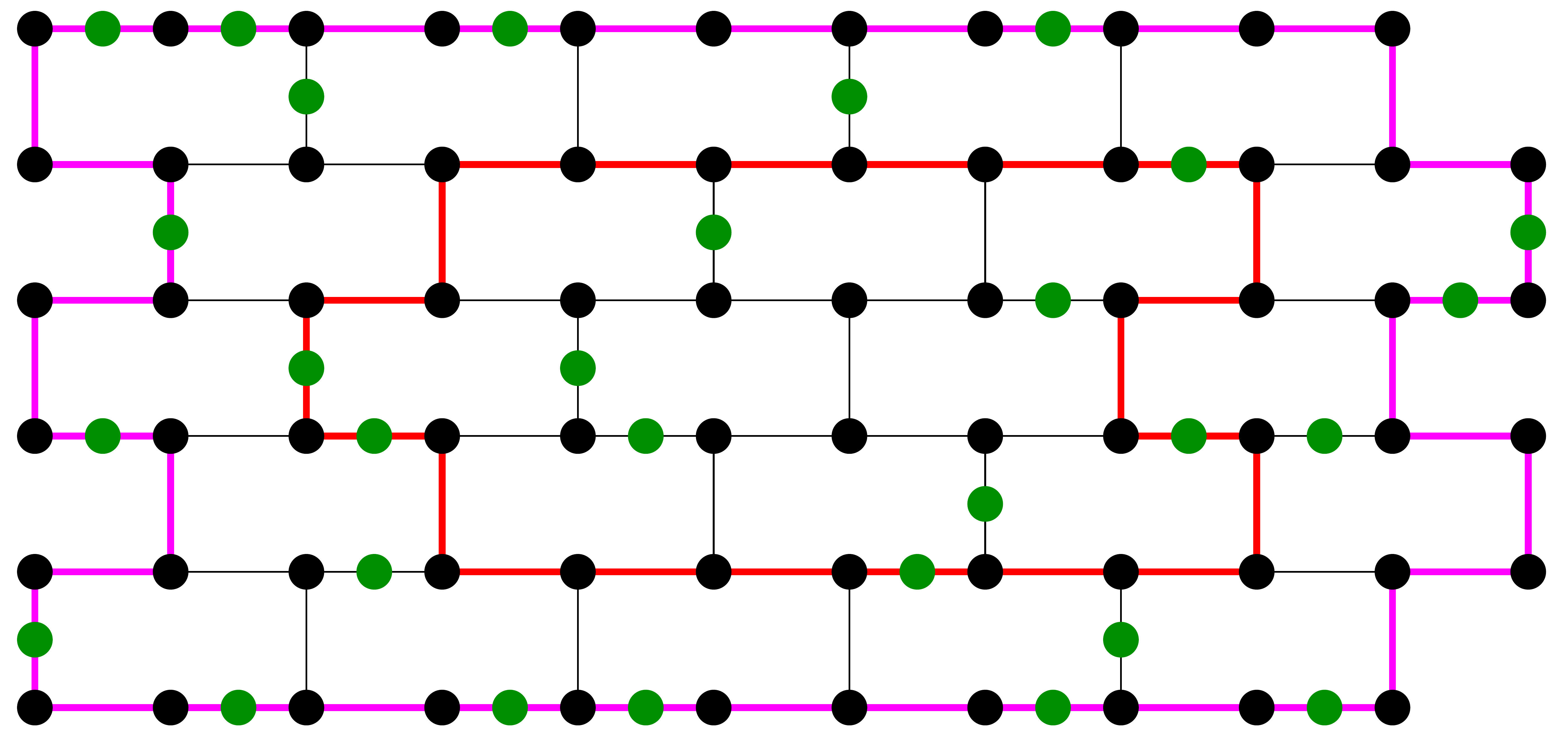}
    \end{center}
%
\caption{The first (magenta) and  second (red) layers of a wall of height 5}
\label{fig:layerw}

\end{figure}

The {\em layers} of a subdivided wall $W$ of height $k$ are recursively defined as follows.
The first layer of $W$ is its perimeter. For $i=2,\cdots,\lfloor \frac{k}{2}\rfloor$, the $i$-th layer of $W$ is the $(i-1)$-th layer of the subwall $W'$ obtained from $W$ after removing from $W$ its perimeter and all occurring vertices of degree 1 (see Figure~\ref{fig:layerw}).

If $W$ is a subdivided wall of height $k$, we call {\em brick} of $W$ any facial cycle whose non-subdivided counterpart in $W_{h}$ has length 6. We say that two bricks are {\em neighbors} if their intersection contains an edge. 

Let $W_{k}$ be a wall. We denote by $P^{(h)}_{j}$ the shortest path connecting vertices $(1,j)$ and $(2\cdot k+2,j)$, and by $P^{(v)}_{i}$ the shortest path connecting vertices $(i,1)$ and $(i,k+1)$ with the assumption that for $i<2\cdot k+2$, $P^{(v)}_{i}$ contains only vertices $(x,y)$ with $x=i,i+1$.
We call the paths $P^{(h)}_{k+1}$ and $P^{(h)}_{1}$ the {\em southern path of $W_{k}$} and {\em northern part of $W_{k}$} respectively. Similarly, we call the paths $P^{(v)}_{1}$ and $P^{(v)}_{2\cdot k+1}$ the {\em western part of $W_{k}$} and the {\em eastern part of $W_{k}$} respectively. If $W$ is a subdivision of $W_{k}$, we will use the same notation for the paths obtained by the subdivisions of the corresponding paths of $W_{k}$.

\paragraph{Compasses and rural devisions.}
Let $W$ be a subdivided wall in $G$.  Let $K'$ be the connected component of $G\setminus P$ that contains $W\setminus P$. The {\em compass} $K$ of $W$ in $G$ is the graph $G[V(K')\cup V(P)]$. Observe that $W$ is a subgraph of $K$ and $K$ is connected. We say that a wall is {\em flat} in $G$ if  its compass $K$ in $G$ has no $(c_{1},c_{3})$-path and $(c_{2},c_{4})$-path
that are vertex-disjoint.

\begin{observation} 
\label{wl:subd:obs}
Let $W$ be a flat wall. Then any subdivision of $W$ is also flat.
\end{observation}

If $J$ is a subgraph of $K$, we denote by $\partial_{K} J$ the set of all vertices $v$ such that either $v\in C$ or $v$ is incident with an edge of $K$ that is not in $J$, i.e. $$\partial_{K} J=\{v\in V(J)\mid v\in C \text{ or } \exists e\in E(K)\setminus E(J): v\in e\}.$$ A {\em rural division} $\mathcal{D}$ of the compass $K$ is a collection $(D_{1},D_{2},\dots,D_{m})$ of subgraphs of $K$ with the following properties:
\begin{enumerate}
\item\label{fldiv1} $\{E(D_{1}),E(D_{2}),\dots,E(D_{m})\}$ is a partition of non-empty subsets of $E(K)$,

\item\label{fldiv2} for $i,j\in [m]$, if $i\neq j$ then $\partial_{K} D_{i}\neq \partial_{K} D_{j}$ and $V(D_{i})\cap V(D_{j})=\partial_{K} D_{i}\cap \partial_{K} D_{j}$,

\item\label{fldiv3} for each $i\in[m]$ and all $x,y\in \partial_{K} D_{i}$ there exists a $(x,y)$-path in $D_{i}$ with no internal vertex in $\partial_{K} D_{i}$,

\item\label{fldiv4} for each $i\in [m]$, $|\partial_{K} D_{i}|\leq 3$, and

\item\label{fldiv5} the hypergraph $\displaystyle H_{K}=(\bigcup_{i\in[m]} \partial_{K} D_{i},\{\partial_{K} D_{i}\mid i\in[m]\})$ can be embedded in a closed disk $\Delta$ such that  $c_{1},c_{2},c_{3}$ and $c_{4}$ appear in this order on the boundary of $\Delta$ and 
for each hyperedge $e$ of $H$ there exist $|e|$ mutually vertex-disjoint paths between $e$ and  $C$ in $K$.
\end{enumerate}

We call the elements of ${\cal D}$ {\em flaps}. A flap $D\in {\cal D}$
is {\em internal} if $V(D)\cap V(P)=\emptyset$.
\\

We are now in position to state the main result of this paper.

\begin{theorem}
\label{th:main}
There exists a computable function $f$ such that for every two graphs $H$ and $G$ and 
every $k\in \mathbb{N}$, one of the following holds:
\begin{enumerate}

\item $H$ is a minor of $G$,

\item $\tw(G)\leq f(h)\cdot k$, where $h=\n(H)$

\item $\exists A\subseteq V(G)$ with $|A|\leq \an(H)-1$ such that $G\setminus A$ contains as a subgraph a flat subdivided wall $W$ where
\begin{itemize} \item  $W$ has height $k$ and
\item the compass of $W$ has  a rural division $\mathcal{D}$ 
such that each internal flap of ${\cal D}$ has  treewidth at most $f(h)\cdot k$.
\end{itemize}
\end{enumerate}
\end{theorem}

\noindent We postpone the proof of Theorem~\ref{th:main} until Section~\ref{sec:proof}.

\section{Some auxiliary lemmata}

The main results in this  section are Lemmata~\ref{wl:subd:inv} and~\ref{apex:ind:claim} that 
will be used for the proof of our main result in Section~\ref{sec:proof}.

\subsection{An invariance lemma for flatness}\label{fil}

\begin{lemma}\label{gd2wl} Let $k$ be a positive integer and $G$ be a graph that contains $\Gamma_{2\cdot k+8}$ as a $v$-smooth contraction. Then $G$ contains as a subgraph a subdivided wall of height $k$ whose compass can be embedded in a closed disk $\Delta$ such that 
the perimeter of $W$ is identical to the boundary of $\Delta$.
\end{lemma}

\begin{proof} Assume that $\Gamma_{2\cdot k+8}$ is a $v$-smooth contraction of $G$ via $\phi$. W.l.o.g. let $V(\Gamma_{2\cdot k+8})=\{1,\dots,2\cdot k+8\}^{2}$, where $v=(2\cdot k+8,2\cdot k+8)$.
Let $R$ be the set of external vertices of $\Gamma_{2\cdot k+8}$ and 
let $G'=G\setminus \bigcup_{x\in R} \phi^{-1}(x)$.
%
As $G$ contains $\Gamma_{2\cdot k+8}$ as a $v$-smooth contraction, it follows that $G'$ is embedded inside an open disk $\Delta'$.
Moreover $G'$ can be contracted to $\Gamma_{2\cdot k+6}^{*}$ via the restriction
of $\phi$ to $V(G')$. From the definition of a wall, it follows that $\Gamma_{2\cdot k+6}^{*}$
contains $W_{k+2}$ as a subgraph.
As $G'$ contains $\Gamma_{2\cdot k+6}^{*}$ as a minor, it follows 
that $G'$ contains $W_{k+2}$ as a minor. As $W_{k+2}$ has maximum degree 3, it is also a topological minor of $G'$. 
Therefore  
$G'$ contains as a subgraph (embedded in $\Delta'$) a subdivided wall of height $k+2$. 
Among all such subdivided walls, let $W_{\rm ex}$ 
be the one whose compass has the minimum number of faces
inside the annulus $\Phi=\Delta_{\rm ex}\setminus \Delta\subseteq \Delta'$ where $\Delta_{\rm ex}$ and $\Delta$ are defined as the closed disks defined so that the boundary of $\Delta_{\rm ex}$ is the first layer of $W_{\rm ex}$ and the boundary of $\Delta$ is the second one.

Let $W$ be the subdivided wall of $G'$ whose perimeter is the boundary of $\Delta$. By definition, all vertices of the compass $K$ of 
$W$ are inside $\Delta$. It now remains to prove that 
the same holds also for the edges of $K$.
Suppose in contrary that $\{x,y\}$ is an edge outside $\Delta$.
Clearly, both $x$ and $y$ lie on the boundary of $\Delta$ and 
$\{x,y\}$ is inside the disk $\Delta^{*}$ defined by some brick of $W_{\rm ex}$ that is inside $\Phi$.  We distinguish two cases:\\

\noindent{\sl Case 1}: $\{x,y\}$ are in the same brick, say $A$ of $W$. Then, there is a path of this brick
that can be replaced in $W$ by $\{x,y\}$ and substitute $W$ by a new subdivided wall corresponding to an annulus with less faces, a contradiction.\\

\noindent{\sl Case 2}: $\{x,y\}$ are not in the same brick of $W$. 
Then $x$ and $y$ should belong in neighboring bricks, say $B$ and $C$ respectively.
Let $A$ be the unique brick of $W_{\rm ex}$ that contains $x$ and $y$
and $w$ be the unique common vertex in $A,B$ and $C$.
Observe that there a path $P_{B}$ of $B$ connecting $x$ and $w$
and a path $P_{C}$ of $C$ connecting $y$ and $w$.
Then we substitute $W$ by a new wall as follows:
we replace $w$ by $x$, $P_{C}$ by $\{x,y\}$, and see $P_{B}$
as a subpath of the common path between $B$ and $C$.
Again, the new wall corresponds to an annulus with less faces, a contradiction.\end{proof}

\begin{lemma}[\cite{FominGT09con}]
\label{remvort}
There is a function $f: \Bbb{N}\times \Bbb{N}\rightarrow \Bbb{N}$ such that if $G$ is a graph $h$-nearly
embedded in a surface of Euler genus $\gamma$ without apices,  where
$\tw(G)\geq f(\gamma,h)\cdot k$, then $G$ contains as a
$v$-smooth contraction the  graph $\Gamma_{k}$
with  the loaded corner $v$.
\end{lemma}

\begin{lemma}
\label{wl:inv}
Let $h$ be a positive integer and $G$ be a graph that contains a flat subdivided wall $W$ of height $h$. If $K_{3}$ is a subgraph of the compass of $W$ then after applying a $\Delta Y$-transformation in $K_{3}$ the resulting graph also contains a flat subdivided wall $W'$ of height $h$ as a subgraph. Moreover, $W'$ is isomorphic to a subdivision of $W$. 
\end{lemma}

\begin{proof} 
We examine the non-trivial case where $E(K_{3})\cap E(W)\neq\emptyset$.
Observe that, as $W$ does not contain triangles, $|E(K_{3})\cap E(W)|<3$. In what follows we denote by $x,y,z$ the vertices of $K_{3}$, $w$ the vertex that appears after the transformation, and distinguish the following cases.

{\sl Case 1.} $K_{3}$ and $W$ have exactly one common edge, say $\{x,y\}$.
As $w$ is a new vertex, the path $(x,w,y)$ that appears after the $\Delta Y$-transformation has no common internal vertices with $W$. In this case, we replace the edge $\{x,y\}$ in $W$ by the edges $\{x,w\},\{w,y\}$.

{\sl Case 2.} $K_{3}$ and $W$ have exactly two common edges, say $\{x,y\}$ and $\{x,z\}$. 
We distinguish the folowing two subcases.

{\sl Subcase 2.1.} $x$ is an original vertex and $x$ is not a corner of $W$. Let $q$ be the third vertex in the neighborhood of $x$. Observe that the $\Delta Y$-transformation is equivalent to removing the edge $\{y,z\}$, which is not an edge of the wall, and subdividing the edge $\{x,q\}$. Then the lemma follows from Observation~\ref{wl:subd:obs}.

{\sl Subcase 2.2.} $x$ is not an original vertex or $x$ is a corner. Then  the  $\Delta Y$-transformation is equivalent to  removing the edge $\{y,z\}$, which is not an edge of $W$, and then substituting $\{y,x\}$ by $\{y,w\}$ and $\{x,z\}$ by $\{w,z\}$. 

Observe that in all the above cases the resulting wall $W'$ remains flat and is isomorphic to a subdivision of $W$ and the lemma follows.
\end{proof}

By applying inductively Observation~\ref{wl:subd:obs} and Lemma~\ref{wl:inv} we derive the following.

\begin{lemma}
\label{wl:subd:inv}
Let $h$ be a positive integer and $G$ be a graph that contains a flat subdivided wall $W$ of height $h$ as a subgraph. If we apply a sequence of 
subdivisions or $\Delta Y$-transformations in $G$,
then the resulting graph will contain a flat subdivided wall $W'$ of height $h$ as a subgraph. Moreover, $W'$ is isomorphic to a subdivision of $W$.
\end{lemma}

\subsection{Pyramids and treewidth}
\label{sec:pytr}

Combining Proposition (1.5) in~\cite{RobertsonST94} with E\"{u}ler's formula for planar graphs we obtain the following.
\begin{lemma}\label{plangm}
If $G$ is a planar graph then  $G$ is isomorphic to a minor of the $(14\cdot \n(G)-24)\times (14\cdot \n(G)-24)$-grid.
\end{lemma}

From the above lemma we obtain the following.
\begin{lemma}\label{apexgraphgrid}
Let $h$ be an integer. If $G$ is an $h$-apex graph then $G$ is isomorphic to a minor of $\Pi_{14\cdot (\n(G) -h)-24,h}$.
\end{lemma}

\begin{lemma}
\label{lem:makepyr}
Let $G$ be the graph obtained by a $((k+\lceil \sqrt{h}\rceil)\times (k+\lceil \sqrt{h}\rceil))$-grid if we make its vertices adjacent to a set $A$ of $h$ new vertices.
Then $G$ contains $\Pi_{k,h}$ as a minor.
\end{lemma}

\begin{proof}
We denote by $G'$ the grid used for constructing $G$ and let $G_{1}$ and 
$G_{2}$ two disjoint subgraphs of $G'$ where $G_{1}$ is isomorphic to 
a $(k\times k)$-grid and $G_{2}$ is isomorphic to a  $(\alpha\times \alpha)$-grid 
where $\alpha=\lceil \sqrt{h}\rceil$. Remove from $G$
all vertices not in $A\cup V(G_{1})\cup V(G_{2})$. 
Then remove all edges of $G'$ incident to $V(G_{2})$
and notice that in the remaining graph $F$, the vertices in $A\cup V(G_{2})$
induce a graph isomorphic to $K_{h,\alpha^{2}}$ which, in turn, can be contracted to a clique on the vertices of $A$. Applying the same contractions in $F$ one may obtain $\Pi_{k,h}$ as a minor of $G$.
\end{proof}

\begin{lemma}
\label{apex:ind:claim}
Ley $G,H$ be graphs such that $H$ is not a minor of $G$ and there exists a set $A\subseteq V(G)$ such that $G\setminus A$ contains a wall $W$ of height $g(h)\cdot (k+1)-1$ as a subgraph, where $g(h)=14\cdot (h-\an(h)) +\lceil \sqrt{\an(h)}\rceil-24$ and $h=\n(H)$.
If $|A|\geq \an(h)$ then there exists an $A'\subseteq A$ such that $|A'|< |A|$ and $G\setminus A'$ contains a wall $W'\subseteq W$ of height $k$ as a subgraph. Moreover, if $K'$ is the compass of $W'$ in $G\setminus A'$ then $V(K')\cap A=\emptyset$. 
\end{lemma}

\begin{proof}
Let $A=\{\alpha_{i}\mid i\in [|A|]\}$ and $P_{g(h)}=\{W_{(m,l)}\mid (m,l)\in [g(h)]^{2}\}$ be a collection of $(g(h))^{2}$ disjoint subwalls $W_{(m,l)},(m,l)\in [g(h)]^{2}$ of $W$ with height $k$.
For every $(m,l)\in [g(h)]^{2}$, we denote by $K_{(m,l)}$ the compass of $W_{(m,l)}$ in $G\setminus A$ and let $q_{(m,l)}=(q_{(m,l)}^{1},q_{(m,l)}^{2},\dots,q_{(m,l)}^{|A|})$ be the binary vector where for every $j\in |A|$, 
$$ q_{(m,l)}^{j}=
\begin{cases} 
1 & \text{if } \exists v\in V(K_{(m,l)}): \{v,\alpha_{j}\}\in E(G)\\
0 & \text{if } \forall v\in V(K_{(m,l)}): \{v,\alpha_{j}\}\notin E(G)
\end{cases}
$$

We claim that there exists an $(m',l')\in [g(h)]^{2}$ such that $q_{(m',l')}\neq (1,1,\dots,1)$. Indeed, assume in contrary, that for every $(m,l)\in [g(h)]^{2}$, $q_{(m,l)}=(1,1,\dots,1)$. We will arrive to a contradiction by showing that $H$ is a minor of $G$. For this, consider the graph $$G'=G[V(W)\cup\bigcup_{(m,l)\in [g(h)]^{2}}V(K_{(m,l)})]\subseteq G.$$ For every $(m,l)\in [g(h)]^{2}$, contract each $K_{(m,l)}$ to a single vertex and this implies  the existence of a  $(g(h)\times g(h))$-grid as a minor of $G'$ and therefore of $G\setminus A$ as well. Moreover,
for each vertex $v$ of this grid it holds that each vertex in $A$ is adjacent to some vertex of the model of $v$, therefore $G$ contains 
the graph $J$ obtained after we take a $(g(h)\times g(h))$-grid
and connect all its vertices with $\an(h)$ new vertices.
From Lemma~\ref{lem:makepyr}, $G$ contains $ \Pi_{14\cdot (n(h)-\an(h))-24, \an(h)}$ as minor. Applying now Lemma~\ref{apexgraphgrid}, we obtain that 
$G$ contains $H$ as a minor, a contradiction.
Therefore, there exist $(m',l')\in [g(h)]^{2}$ and $j_{0}\in [|A|]$ such that $q_{(m',l')}^{j_{0}}=0$. The lemma follows for $A'=A\setminus \{\alpha_{j_{0}}\}$ and $W'=W_{(m',l')}$.
\end{proof}

\section{The main proof}

\label{sec:proof}

Given a tree decomposition ${\cal T}=({\cal X}=\{X_{i}\mid i\in V(T)\},T)$ of a graph $G$ a vertex $i_{0}\in V(T)$ 
and a set of vertices ${I}\subseteq N_{T}(i_{0})$, 
we define ${\cal T}_{i_{0},{ I}}$ as the collection of connected 
components of $T\setminus  i_{0}$ that contain vertices of ${ I}$.
Given a subtree $Y$ of $T$, we define $G_{Y}=G[\cup_{i\in V(Y)}X_{i}]$
and $\overline{G}_{Y}=\cup_{i\in V(Y)}\overline{X}_{i}$. 

\begin{observation}
\label{obs:complete}
Given a tree decomposition ${\cal T}=({\cal X}=\{X_{i}\mid i\in V(T)\},T)$ of a graph $G$, a vertex $i_{0}\in V(T)$, and a set of vertices ${I}\subseteq N_{T}(i_{0})$, it holds that  
\mbox{for every ${T}_{1},{T}_{2}\in {\cal T}_{i_{0},{I}}$},  $\overline{G}_{{T}_{1}}\cap \overline{G}_{{T}_{2}}$ \mbox{is a complete graph}.
\end{observation}

Given a family of graphs ${\cal F}$, a graph $G$ and a set of vertices $S\subseteq V(G)$, we  define the class ${\cal F}_{S,G}^*$ as the collection of the connected components in the graphs of  ${\cal F}\setminus S$ and the class ${\cal F}_{S,G}$ as the set of graphs in ${\cal F}_{S,G}^*$ that have some common vertex with $G\setminus S$. We say that two graphs $G_{1},G_{2}\in {\cal F}_{S,G}$ are $G$-equivalent if $V(G_{1})\cap V(G\setminus S)=V(G_{2})\cap V(G\setminus S)$ and let ${\cal F}_{S,G}^{1},\dots,{\cal F}_{S,G}^{\rho}$ be the equivalence classes defined that way. We denote by ${\cal P}_{{\cal F},S,G}=\{\cupall {\cal F}_{S,G}^{1},\dots,\cupall {\cal F}_{S,G}^{\rho}\}$,  i.e. for each equivalence class ${\cal F}_{S,G}^{i}$ we construct a graph in ${\cal P}_{{\cal F},S,G}$, by taking  the union of the graphs in ${\cal F}_{S,G}^{i}$.\\

\begin{proof}[Proof of Theorem~\ref{th:main}]
Let $G$ be a graph that excludes $H$ as a minor. By Proposition~\ref{structurethm}, there is a computable function $f_{1}$ such that 
there exists a tree decomposition $${\cal T}=({\cal X}=\{X_{i}\mid i\in V(T)\},T)$$ of $G$, where  for every $i\in V(T)$, the graphs $\overline{X}_{i}$ are $f_{1}(h)$-nearly-embeddable in a surface $\Sigma$ of genus $f_{1}(h)$. 
Among all such tree-decompositions we choose ${\cal T}=({\cal X},T)$ such that 
\begin{enumerate}[(i)]
\item\label{td:def1} ${\cal T}$ is small.
\item\label{td:def2}  subject to (i) ${\cal T}$ has maximum number of nodes.
\item\label{td:def3} subject to (ii) the quantity $\displaystyle 
\sum_{\substack{i,j\in V(T)\\ i\neq j}}|X_{i}\cap X_{j}|$ is minimized.
\end{enumerate}
Notice that, from Lemma~\ref{fincomptd},  the condition~(\ref{td:def1}) guaranties the possibility of the choice of condition (ii). We use the notation $\overline{G}$ to denote the  graph $\overline{G}_{T}$ and we call the edges of $E(\overline{G})\setminus E(G)$ {\em virtual}.


Let $w:V(T)\rightarrow \mathbb{N}$ such that $w(i)=\tw(\overline{X}_{i})$. Observation~\ref{verwithungrtw} and Observation~\ref{treew:bags} imply that there exists a vertex $i_{0}\in V(T)$ such that $\tw(\overline{X}_{i_{0}})\geq \tw(G)$ and at most one of the connected components of $T\setminus i_{0}$ contains a vertex $j$ such that $w(j)> w(i_{0})$. We denote by ${A}_{i_{0}}$ the set of apices of the graph $\overline{X}_{i_{0}}$ and by $F$ the graph $\overline{X}_{i_{0}}\setminus {A}_{i_{0}}$ (notice that $F\subseteq \overline{G}$ but $F$ is not necessarily a subgraph of $G$ as $F$ may contain virtual edges).\medskip

\noindent{\sl Claim 1.} For every connected component $Y$ of $T\setminus \{i_{0}\}$, there is a vertex in $G_{Y}\setminus X_{i_{0}}$ connected with $X_{i_{0}}\cap V(G_{Y})$
with $|X_{i_{0}}\cap V(G_{Y})|$ vertex disjoint paths whose internal vertices belong in $G_{Y}\setminus X_{i_{0}}$.
 
 \noindent{\sl Proof of Claim 1.} First, observe that, from~(\ref{td:def2}) the graph $G_{Y}\setminus X_{i_{0}}$ is connected. Then, from~(\ref{td:def3}), it follows that there is a vertex in $G_{Y}\setminus X_{i_{0}}$ connected with $X_{i_{0}}\cap V(G_{Y})$ with $|X_{i_{0}}\cap V(G_{Y})|$ vertex disjoint paths.\medskip

From Lemma~\ref{apex} and the choice of $i_{0}$ holds that 
\begin{equation}\label{ineq1}\tw(F)=\tw(\overline{X}_{i_{0}}\setminus {A}_{i_{0}})\geq\tw(\overline{X}_{i_{0}})- |A_{i_{0}}|\geq \tw(G) - |A_{i_{0}}|.\end{equation}
 Recall that $|A_{i_{0}}|\leq f_{1}(h)$.  
Let $f_{2}$ be the two-variable function of Lemma~\ref{remvort}.
We define  the two-variable function $f_{3}$  and the one-variable functions $f_{4}$ and $f_{5}$ such that
\begin{eqnarray*}
f_{5}(h) & = & 14\cdot (h-\an(h)) +\lceil \sqrt{\an(h)}\rceil-24\\
f_{4}(h) & = & f_{5}(h)^{|A_{i_{0}}|-\an(h)+1}\\
f_{3}(h,k) & = & f_{2}(f_{1}(h),f_{1}(h))\cdot (4k\cdot f_{4}(h)+12)+f_{1}(h)
\end{eqnarray*}

As $F$ is $f_{1}(h)$-nearly embeddable in $\Sigma$ and does not contain any apices, from \eqref{ineq1} and Lemma~\ref{remvort}, we obtain  that  if $\tw(G)\geq f_{3}(h,k)$ then $F$ contains the graph $Q=\Gamma_{4k\cdot f_{4}(h)+12}$ as a $v$-smooth contraction.
By Lemma~\ref{gd2wl},  it follows that $F$ contains as a subgraph a flat subdivided wall $W'$ of height $2k\cdot f_{4}(h)+2$ whose compass $K'$ in $F$ can be embedded  in a closed disk $\Delta$ such that 
the perimeter of $W'$ is identical to the boundary of $\Delta$. Let
 $${I}'=\{i\in N_{T}(i_{0})\mid X_{i}\cap V(K')\neq \emptyset\}.$$ 
In other words $I'$ corresponds to all nodes of the tree decomposition ${\cal T}$ that have vertices ``inside'' the compass of the subdivided wall $W'$ in $F$.
Clearly, every clique in $K'$ has size at most $4$. Furthermore, $K'$ does not contain any $K_{4}$ as a subgraph. Indeed,
if so, one of its triangles would be a separator of $G$, contradicting minimality of
condition (iii) in the definition of the decomposition ${\cal T}$. To see this,  just replace $X_{i_{0}}$ in ${\cal T}$ by $X_{i_{0}}\setminus \{z\}$ where $z$ is the vertex of the clique that is not in the separating triangle. Recall, now, that for every $i\in \mathcal{I}'$, $F[V(K')\cap X_{i}]\subseteq K'$ is a clique. Thus,
 \begin{equation}\label{eqncli} \text{for every }i\in I',|V(K')\cap X_{i}|\leq 3.\end{equation}
 %
%
Recall, also, that there exists at most one tree in ${\cal T}_{i_{0},{I}'}$, say $T'$, that contains a vertex $i_{1}$ with $w(i_{1})> w(i_{0})$. 
%
%
Let ${\cal W}'=\{W_{1}', W_{2}', W_{3}', W_{4}'\}$ be the collection of 
vertex disjoint subwalls of $W'$ of height $f_{4}(h)\cdot k$ not meeting the vertices of $P_{k\cdot f_{4}(h)+2}^{(h)}$ and $P_{2k\cdot f_{4}(h)+3}^{(v)}$ (see Figure~\ref{fig:pthw}).

\begin{figure}[h]
  \begin{center}
\includegraphics[width=11cm]{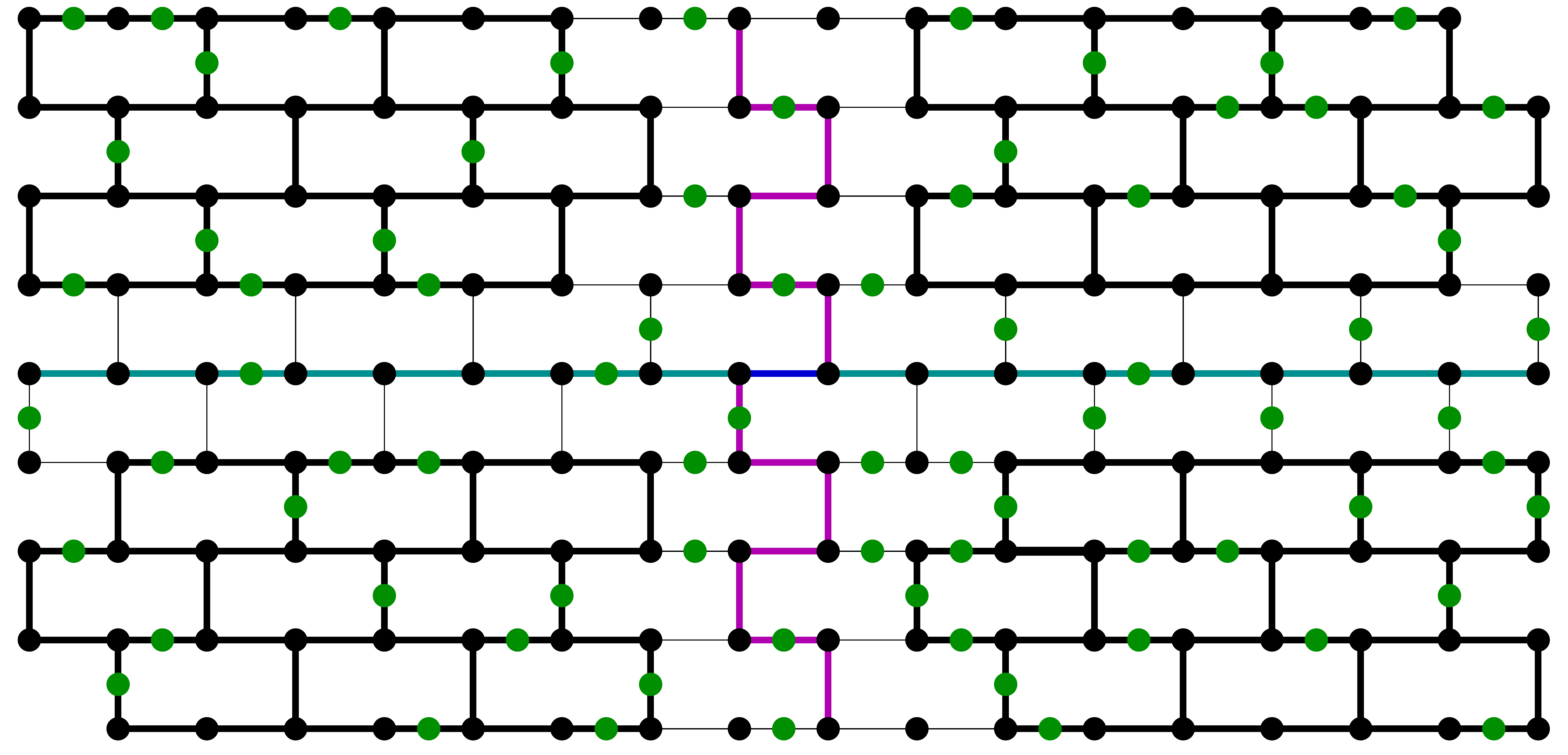}
    \end{center}
\label{fig:pthw}
\caption{The paths $P_{k\cdot f_{4}(h)+2}^{(h)}$ (cyan) and $P_{2k\cdot f_{4}(h)+3}^{(v)}$ (magenta) and the corresponding walls for $k=1$ and $f_{4}(h)=3$.}
\end{figure}
 
From \eqref{eqncli}, $X_{i_{1}}$ has at most 3 vertices in common with $K'$, therefore there exists a subwall $\tilde{W}\in{\cal W}'$  of height $f_{4}(h) \cdot k$, with compass $\tilde{K}$ in $F$ such that $V(\tilde{K})\cap V(G_{T'})=\emptyset$. Consequently, if we set 
$$\tilde{I}=\{i\in N_{T}(i_{0})\mid X_{i}\cap V(\tilde{K})\neq \emptyset\}$$
we have that $\tilde{I}\subseteq I'\setminus \{i_{1}\}$ and for every tree $\tilde{T}\in {\cal T}_{i_{0},\tilde{I}}\subseteq {\cal T}_{i_{0},I'}\setminus\{T'\}$ it holds 
that $\max\{w(i)\mid i\in V(\tilde{T})\}\leq f_{3}(h,k)$. Therefore, for every $\tilde{T}\in {\cal T}_{i_{0},\tilde{I}}$, $\tw(\overline{G}_{\tilde{T}})\leq f_{3}(h,k)$.
As $G_{\tilde{T}}$ is a subgraph of $\overline{G}_{\tilde{T}}$, it follows that 
\begin{eqnarray}
\mbox{for every $\tilde{T}\in {\cal T}_{i_{0},\tilde{I}}$,} & \tw(G_{\tilde{T}})\leq f_{3}(h,k).
\label{eq:twmo}\end{eqnarray}

From~\eqref{eqncli}, it follows that for every $\tilde{T}\in {\cal T}_{i_{0},\tilde{I}}$,
 the vertices in $V(G_{\tilde{T}})\cap V(\tilde{K})$ induce a clique in $\tilde{K}$  with at most 3 vertices where some of its  edges may be virtual. \\

Let $\tilde{V}=V(F)\setminus V(\tilde{K})$ and ${\cal F}'=\{G_{\tilde{T}}\mid \tilde{T}\in{\cal T}_{i_{0},\tilde{I}}\}$. Notice that $\tilde{K}=F\setminus \tilde{V}$.
We denote by ${\cal F}$ the class ${\cal P}_{{\cal F}',\tilde{V},F}$.

We call the edges in  $\tilde{E}=E(\tilde{K})\setminus E(G)$ {\em useless}.  We also call all vertices 
in $V(\cupall {\cal F})\setminus V(\tilde{K})$ {\em flying} vertices.
The non-flying vertices of a graph $R$ in ${\cal F}$ are the {\em base} of $R$. 
Notice that, by  the definition of ${\cal F}$, each graph $R$ in ${\cal F}$ 
is a subgraph of the union of some graphs of ${\cal F}'$.
From~Observation~\ref{obs:complete} and  \eqref{eq:twmo},
It follows that 
\begin{itemize}
\item[(a)] all graphs in ${\cal F}$ have treewidth at most $f_{3}(h,k)$
\end{itemize}
Observation~\ref{obs:complete} and  \eqref{eqncli} yields that
\begin{itemize}
\item[(b)] the base vertices of each $R$ induce a clique of size 1,2, or 3 in $\tilde{K}$.
\end{itemize}
Also, from Claim 1 and the fact that $\tilde{V}\cup A_{i_{0}}\subseteq X_{i_{0}}$, we have that
\begin{itemize}
\item[(c)] each pair of vertices of some graph in ${\cal F}$
are connected in $G$ by a path whose internal vertices are flying.
\end{itemize}
Note that each clique mentioned  in (b) may contain useless edges. Moreover, from (c), all virtual edges of $\tilde{K}$ are edges of such a clique.
Let $\tilde{G}=(V(G),\tilde{E}\cup E(G))$, i.e. we 
add in $G$ all useless edges.

It now follows that $\tilde{G}\setminus A_{i_{0}}$ contains the wall $\tilde{W}$ as a subgraph and the compass 
of $\tilde{W}$ in $\tilde{G}\setminus A_{i_{0}}$ is 
$$\tilde{K}^{+}=\tilde{K}\cup\bigcup_{R\in{\cal F}}R $$

Notice that the wall $\tilde{W}$ remains flat in $\tilde{G}$. Indeed, suppose that $Q_{1}$ and $Q_{2}$ are two vertex disjoint paths between the two anti-diametrical corners of $\tilde{W}$ such that the sum of their lengths is minimal. As not both of $Q_{1}$ and $Q_{2}$ may exist in $\tilde{K}$, some of them, say $Q_{1}$ contains some flying vertex. Let $R$ be the graph in ${\cal F}$ containing that vertex. Then there are two vertices $x$ and $y$ of the base of $R$ 
met by $Q_{1}$.  From (b), $\{x,y\}$ is an edge of $\tilde{K}^{+}$ and we can substitute the portion of $Q_{1}$ that contains flying vertices by $\{x,y\}$, a contradiction to the minimality of the choice of $Q_{1}$ and $Q_{2}$. \\

Let $\tilde{E}^{+}=E(\tilde{K}^{+})\setminus E(\cupall {\cal F})$, i.e. $\tilde{E}^{+}$ is the set of edges of $\tilde{K}$ not contained in any graph $R$ of ${\cal F}$. It follows that all useless edges are contained in $\tilde{E}^{+}$, i.e.
\begin{equation}
\tilde{E} \subseteq \tilde{E}^{+}\label{eq6}
\end{equation}
For every $e\in \tilde{E}^{+}$, we denote by $\tilde{G}_{e}$ the graph formed by the edge  $e$ (i.e. the graph $\tilde{G}[e]$) and let ${\cal E}=\{\tilde{G}_{e}\mid e\in \tilde{E}^{+}\}$. We set $\tilde{{\cal D}}^{+}={\cal F}\cup {\cal E}$. Notice that, 
\begin{eqnarray}
\text{For every graph } R\in{\cal F}, & & \partial_{\tilde{K}^{+}} R \text{ is the base of } R \label{eq4}\\
\text{For every graph } \tilde{G}_{e}\in {\cal E}, & & \partial_{\tilde{K}^{+}} \tilde{G}_{e}=V(\tilde{G}_{e}) \label{eq5}
\end{eqnarray}
\medskip

\noindent{\sl Claim 2}. $\tilde{{\cal D}}^{+}={\cal F}\cup {\cal E}$ is a rural division of $\tilde{K}^{+}$. 

\noindent{\sl Proof of Claim 2}. Properties~\ref{fldiv1} and~\ref{fldiv2}, follow from the construction of the graphs in ${\cal F}$ and ${\cal E}$. Moreover, Properties~\ref{fldiv3} and~\ref{fldiv4}
follow from  (c)  and (b) respectively. 
For Property~\ref{fldiv5}, recall that $\tilde{W}$ is a subwall of $W'$ whose compass $K'$ in $F$ can be embedded in a closed disk $\Delta$ such that the perimeter of $W'$ is identical to its boundary. This implies that $\tilde{K}$ can be embedded in a closed disk $\tilde{\Delta}\subseteq \Delta$ such that the corners $c_{1},c_{2},c_{3}$, and $c_{4}$ of $\tilde{W}$ appear in this order on its boundary. We now consider the following hypergraph:
 $$\tilde{H}^{+}=(\cupall\{\partial_{\tilde{K}^{+}}D\mid D\in \tilde{{\cal D}}^{+}\},\{\partial_{\tilde{K}^{+}} D\mid D \in \tilde{{\cal D}}^{+}\}).$$  
%
Notice that $V(\tilde{H}^{+})=V(\tilde{K})$. We can now construct $I(\tilde{H}^{+})$ by applying, for each $D\in \tilde {\cal D}^{+}$, the following transformations on the planar graph $\tilde{K}$.

\begin{itemize}
\item  If $|\partial_{\tilde{K}^{+}} D|=1$, we add a new vertex  and an edge that connects it the unique  vertex of $\partial_{\tilde{K}^{+}} D$.

\item  If $|\partial_{\tilde{K}^{+}} D|=2$,  we subdivide the edge of $\tilde{K}[\partial_{\tilde{K}^{+}} D]$ (recall that $\tilde{K}[\partial_{\tilde{K}^{+}}D]$ is isomorphic to $K_{2}$).

\item  If $|\partial_{\tilde{K}^{+}} D|=3$, we apply a $\Delta Y$-transformation in $\tilde{K}[\partial_{\tilde{K}^{+}}D]$ (recall that $\tilde{K}[\partial_{\tilde{K}^{+}} D]$ is isomorphic to $K_{3}$).
\end{itemize}

From Observation~\ref{wl:subd:obs} and Lemma~\ref{wl:inv} follows that the obtained graph remains embedded in $\tilde{\Delta}$. 
It now remains to show that for each $e\in E(\tilde{H}^{+})$ there exist $|e|$ vertex-disjoint paths between $|e|$ and $C$ in $\tilde{K}^{+}$.
Notice that for each $e\in E(H^{+})$ the vertices of $e$ belong to $\tilde{K}$. Finally, there
are  $|e|$ paths between $e$ and $C$, otherwise we would have a contradiction to the choice of the tree-decomposition (assumption (iii)). Therefore all conditions required for Claim 2 hold.
\medskip

Our aim now is to find in $G\setminus A_{i_{0}}$ a flat subdivided wall $\widehat{W}$ of height $f_{4}(h)\cdot k$. From (b),(c), and~\eqref{eq4},  all the useless edges of $\tilde{K}$ are induced by the sets $\partial_{\tilde{K}^{+}} R$, $R \in {\cal F}$ where $\tilde{K}[\partial_{\tilde{K}^{+}} R]$ is isomorphic to either $K_{2}$ or $K_{3}$. Our next step is to prove that, in both cases, we may find a flat subdivided wall in $G\setminus A_{i_{0}}$ of height $f_{4}(h)\cdot k$ that does not contain any useless edges.

{\em Case 1}.~$\tilde{K}[\partial_{\tilde{K}^{+}} R]$ is isomorphic to $K_{2}$.~Then, from (c), there exists a path in $R$ whose endpoints are the vertices of $\partial_{\tilde{K}^{+}} R$ and such that its internal vertices are flying. 

{\em Case 2}.~$\tilde{K}[\partial_{\tilde{K}^{+}} R]$ is isomorphic to $K_{3}$.~Claim 1, combined with the facts that $\tilde{V}\cup A_{i_{0}}\subseteq X_{i_{0}}$ and that $\forall_{R\in {\cal F}}\ \partial_{\tilde{K}^{+}} R\subseteq X_{i_{0}}$, imply that there exists a flying vertex 
$v_{R}$ in $R$ and vertex-disjoint paths between $v_{R}$ and the vertices of $\partial_{\tilde{K}^{+}} R$ whose internal vertices are also flying.

The above case analysis implies that for each $R\in {\cal F}$ the edge $\{x,y\}$ or the triangle with vertices 
$\{x,y,z\}$, induced by 
$\partial_{\tilde{K}^{+}} R$ may be substituted, 
using subdivisions or $\Delta Y$-transformations by a flying path between $x$ and $y$ or by three flying paths 
from a flying vertex $v_{R}$ to $x,y,$ and $z$ respectively. As all edges of these paths 
are flying, they cannot be useless and therefore they exist also in $G\setminus A_{i_{0}}$.
We are now in position to apply 
Observation~\ref{wl:subd:obs} and Lemma~\ref{wl:subd:inv} and obtain  that  $\tilde{G}\setminus A_{i_{0}}$ contains  a flat subdivided wall $\widehat{W}$ of height $f_{4}(h)\cdot k$ such that
\begin{enumerate}[(I)]
\item $E(\widehat{W})\cap \tilde{E}=\emptyset$ (recall that $\tilde{E}$ is the set of the useless edges) and
\item $\widehat{W}$ is isomorphic to a subdivision of $\tilde{W}$.
\end{enumerate}
Therefore, from (I), $\widehat{W}$ is a flat subdivided wall of height $f_{4}(h)\cdot k$ in $G\setminus A_{i_{0}}$. \medskip

Let $\tilde{C}$ and $\widehat{C}$ be the corners of $\tilde{W}$ and $\widehat{W}$ respectively.
We denote by $\sigma$ be the bijection from $\tilde{C}$ to $\widehat{C}$ induced by the 
isomorphism in (II).  
We also enhance $\sigma$ by defining $\phi=\sigma\cup\{(x,x)\mid x\in V(\tilde{W})\setminus C(\tilde{W})\}$. 

Let $\widehat{K}$ be the compass of $\widehat{W}$ in $G\setminus A_{i_{0}}$. 
We claim that

$$\widehat{{\cal D}}=\{D\cap \widehat{K}\mid  D\in \tilde{\cal D}^{+}\}$$
is a 
%
 rural division of $\widehat{K}$. This is easy to verify in what concerns conditions (1--4).
Condition (5) follows by the observation that the mapping $\phi$, defined above, is an isomorphism between $H_{\tilde{K}^{+}}$ and $H_{\widehat{K}}$.

So far, we have found a flat subdivided wall $\widehat{W}$ in $G\setminus A_{i_{0}}$ and a rural division of its compass $\widehat{K}$. As each flap in $\widehat{\cal D}$ is a subgraph of a flap in $\tilde{{\cal D}^{+}}$ we obtain that all flaps in $\widehat{\cal D}$ have  treewidth at most $f_{3}(h,k)$. 
By applying Lemma~\ref{apex:ind:claim} $|A_{i_{0}}|-\an(h)+1$ times, it follows that there exists a set $A\subseteq A_{i_{0}}$, such that $|A|\leq \an(h)-1$ and $G\setminus A$ contains a flat subdivided wall $W$ of height $k$ such that $W\subseteq \widehat{W}$. Moreover, $V(K)\cap A_{i_{0}}=\emptyset$, where $K$ is the compass of $W$ in $G\setminus A$. As above,
$${\cal D}=\{D\cap K\mid D\in \widehat{\cal D}\}$$
is a rural devision of $K$ where all of its flaps have treewidth at most $ f_{3}(h,k).$ The theorem follows as 
 $f_{3}$ is a linear function of $k$. 
\end{proof}
\medskip

The following corollary gives a more precise  description of the structure of apex minor free graphs.\\

\begin{corollary}
There exists a computable function $f$ such that for every two graphs $H$ and $G$, where $H$ is an apex graph and 
every $k\in \mathbb{N}$, one of the following holds:
\begin{enumerate}
\item $\tw(G)\leq f(h)\cdot k$, where $h=|V(H)|$

\item $H$ is a minor of $G$,

\item $G$ contains a flat subdivided wall $W$ where 
\begin{itemize} 
\item $W$ has height $k$ and
\item the compass of $W$ has a rural division $\mathcal{D}$ such that each internal flap of $\mathcal{D}$ has treewidth at most $ f(h)\cdot k$.
\end{itemize}
\end{enumerate}
\end{corollary}
%
%

%
\end{document}